\documentclass[reqno]{amsart}
\usepackage{amsmath,amsthm,amssymb,amscd,float}
\usepackage{subcaption}
\usepackage[shortlabels]{enumitem}
\usepackage{youngtab}
\usepackage{tikz}
\usepackage{xcolor}
\usetikzlibrary{arrows,matrix}
\usepackage{amsfonts,hyperref} 
\usepackage{tabularx}
\newtheorem{theorem}{Theorem}[section]

\newtheorem{conjecture}[theorem]{Conjecture}
\newtheorem{remark}{Remark}[section]
\usepackage[english]{babel}
\usepackage{dutchcal}
\allowdisplaybreaks

\usepackage{ulem}

\title[On hook length biases in $t$-regular partitions]{On hook length biases in $t$-regular partitions}

\author[R. Barman]{Rupam Barman}
\address{Department of Mathematics, Indian Institute of Technology Guwahati, Assam, India, PIN- 781039}
\email{rupam@iitg.ac.in}

\author[P. J. Mahanta]{Pankaj Jyoti Mahanta}
\address{Gonit Sora, Dhalpur, Assam 784165, India}
\email{pjm2099@gmail.com}

\author[G. Singh]{Gurinder Singh}
\address{Department of Mathematics, Indian Institute of Technology Guwahati, Assam, India, PIN- 781039}
\email{gurinder.singh@iitg.ac.in}



\begin{document}

\begin{abstract}
Let $t\geq2$ and $k\geq1$ be integers. A $t$-regular partition of a positive integer $n$ is a partition of $n$ such that none of its parts is divisible by $t$. Let $b_{t,k}(n)$ denote the number of hooks of length $k$ in all the $t$-regular partitions of $n$. Recently, the first and the third authors proved that $b_{3,2}(n)\geq b_{2,2}(n)$ for all $n\geq 4$, and conjectured that $b_{t+1,2}(n)\geq b_{t,2}(n)$ for all $t\geq 3$ and $n\geq 0$. In this paper, we prove that the conjecture is true for $t=3$.
\end{abstract}

\maketitle

\section{Introduction}
A partition $\lambda$ of a positive integer $n$ is a finite sequence of positive integers $\lambda=(\lambda_1, \lambda_2, \ldots, \lambda_r)$ such that $\lambda_1\geq \lambda_2\geq \cdots \geq \lambda_r$ and $\sum\limits_{i=1}^r\lambda_i = n$. The numbers $\lambda_1, \lambda_2, \ldots, \lambda_r$ are called the parts of the partition $\lambda$. Let $t\geq 2$ be a fixed positive integer. A $t$-regular partition of a positive integer $n$ is a partition of $n$ none of whose parts is divisible by $t$.
\par
\begin{figure}[h]
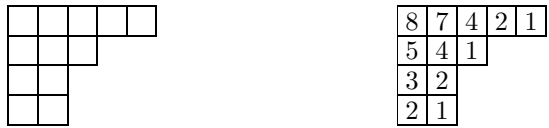

\centering
\begin{minipage}[b]{0.4\textwidth}
\[\young(~~~~~,~~~,~~,~~)\]
\end{minipage}
\begin{minipage}[b]{0.4\textwidth}
\[\young(87421,541,32,21)\]
\end{minipage}
\caption{The Young diagram of the partition $(5,3,2,2)$ and its hook lengths}\label{Figure6.1}
\end{figure}

 
A \textit{Young diagram} of a partition $(\lambda_1, \lambda_2, \ldots, \lambda_r)$ is a left-justified array of boxes, where the $i$-th row from the top contains $\lambda_i$ boxes. For example, the Young diagram of the partition $(5,3,2,2)$ is shown in Figure \ref{Figure6.1} (left). The \textit{hook length} of a box in a Young diagram is the sum of the number of the boxes directly right to it, the number of boxes directly  below it and 1 (for the box itself). For example, see Figure \ref{Figure6.1} (right) for the hook lengths of each box in the Young diagram of the partition $(5,3,2,2)$. A hook of length $k$ is also called $k$-hook.
 
\par 
In recent times, several authors have studied hook lengths of various partition functions. Other than the ordinary partition function, there have been a lot of studies on hook lengths of several restricted partition functions, for example, partitions into odd parts, partitions into distinct parts, partitions into odd and distinct parts, self conjugate partitions and doubled distinct partitions, see e.g. \cite{Ono_Singh, Ballantine_2023,Craig,Han_2016,Han_2017,Petreolle,Singh1Barman}. First of its kind, Ballantine et al. \cite{Ballantine_2023} proved certain hook length biases between partitions into odd parts and partitions into distinct parts. Following Ballantine et al., several authors have worked on hook length biases, see e.g., \cite{Cossaboom,Craig,Singh1Barman,Singh2Barman}.
\par For integers $t\geq2$ and $k\geq1$, let $b_{t,k}(n)$ denote the number of hooks of length $k$ in all the $t$-regular partitions of $n$. 
Recently, the first and the third authors \cite{Singh1Barman, Singh2Barman} studied hook length biases among $b_{t,k}(n)$ for certain values of $t$ and $k$. In \cite{Singh1Barman}, the authors studied the hook length biases for $2$- and $3$-regular partitions for different hook lengths. They established two hook length biases for $2$-regular partitions, namely, $b_{2,2}(n)\geq b_{2,1}(n)$, for all $n>4$ and $b_{2,2}(n)\geq b_{2,3}(n)$, for all $n\geq0$. In \cite{Singh2Barman}, the authors studied biases among $b_{t,k}(n)$ for fixed $k$. In one of their main results, they proved that $b_{3,2}(n)\geq b_{2,2}(n)$ for all $n\geq 4$. Based on numerical evidences, they further proposed the following conjecture.
\begin{conjecture}[Conjecture 4.1, \cite{Singh2Barman}]\label{con1}
Let $t\geq 3$ be an integer. We have $b_{t+1,2}(n)\geq b_{t,2}(n)$, for all $n\geq 0$.
\end{conjecture}
In this article, we prove the following theorem, which confirms the truth of Conjecture \ref{con1} for $t=3$.
\begin{theorem}\label{Th1}
Conjecture \ref{con1} is true for $t=3$.
\end{theorem}
The authors in \cite{Ballantine_2023,Cossaboom,Craig,Singh1Barman,Singh2Barman} proved the hook length biases using combinatorial techniques and the circle method. Our method of proof of Theorem \ref{Th1} is completely combinatorial. 
\section{Proof of Theorem \ref{Th1}}\label{sec:proof}
First, we recall another form of representation of a partition given by $$\lambda=(\lambda^{m_1}_1,\lambda^{m_2}_2,\ldots,\lambda^{m_r}_r),$$ where $m_i$ is the multiplicity of the part $\lambda_i$ and $\lambda_1>\lambda_2>\cdots>\lambda_r$.
Let $\mathbcal{B}_t(n)$ denote the set of all the $t$-regular partitions of $n$. We arrange the parts of a partition $\tau$ in $\mathbcal{B}_3(n)$ into blocks as follows:
\begin{align*}
	\tau= &\bigg((12k+11)^{\alpha_{k,11}}, && (12k+10)^{\alpha_{k,10}},   &&&&\leftarrow 3_k\text{-block}\\
	&\hspace{3mm} (12k+8)^{\alpha_{k,8}},&& (12k+7)^{\alpha_{k,7}}, && (12k+5)^{\alpha_{k,5}},  &&\leftarrow 2_k\text{-block}\\
	&\hspace{3mm} (12k+4)^{\alpha_{k,4}}, && (12k+2)^{\alpha_{k,2}}, &&  (12k+1)^{\alpha_{k,1}}\bigg)_{k\geq 0},  &&\leftarrow 1_k\text{-block}
\end{align*}
where $\alpha_{k,j}$ denotes the multiplicity of the part $12k+j$.
\par 
We consider the map $\Phi:\mathbcal{B}_3(n)\to\mathbcal{B}_4(n)$ defined by
\begin{align*}
		\Phi(3_k\text{-block}) & = \left((12k+11)^{\alpha_{k,11}},\quad (12k+10)^{\alpha_{k,10}}\right),\\
	\Phi(2_k\text{-block}) & = \left((12k+7)^{\alpha_{k,7}},  \quad (12k+6)^{\alpha_{k,8}}, \quad (12k+5)^{\alpha_{k,5}}\right),
\end{align*}
\begin{equation*}
    \Phi(1_k\text{-block}) =\begin{cases}
	   \left((12k+3)^{\alpha_{k,4}},\quad (12k+2)^{\alpha_{k,2}},\quad (12k+1)^{\alpha_{k,1}}\right)\quad &\text{ if } k>0,\\
		\left(3^{\alpha_{0,4}},\quad 2^{\alpha_{0,2}},\quad 1^{\alpha_{0,1}+\sum_{j\geq 0}\left(2\alpha_{j,8} + \alpha_{j,4}\right)}\right)&\text{ if } k=0.
	\end{cases} 
\end{equation*}
For example, $\Phi((10^2, 8^3, 5, 4, 2, 1^2))=(10^2, 6^3, 5, 3, 2, 1^9)$.
The map $\Phi$ is injective. It is clear that we can apply the map $\Phi$ to the partitions in $\mathbcal{B}_3(n)$ for any $n$.
\par 
Let $\mathbcal{h}(i_k\text{-block})$ and $\mathbcal{h}(\Phi(i_k\text{-block}))$ denote the number of 2-hooks in the $i_k\text{-block}$ and $\Phi(i_k\text{-block})$ respectively for $k\geq 0$ and $i\in\{1,2,3\}$. Now, we study the loss of 2-hooks in the partitions under the map $\Phi$. We consider the following three cases.\\
\noindent\textbf{Case 1: $\sum_{k\geq 0}\left(\alpha_{k,8} + \alpha_{k,4}\right)=0$.}
That is, $\alpha_{k,8}=0$ and $\alpha_{k,4}=0$ for all $k\geq 0$. In this case, $\Phi(\tau)=\tau$, and there is no loss of 2-hooks under $\Phi$.\\
\noindent\textbf{Case 2: $\sum_{k\geq 0}\left(\alpha_{k,8} + \alpha_{k,4}\right)\neq0$ and $\sum_{k\geq 0}\left(\alpha_{k,5} + \alpha_{k,2}\right)=0$.}
In this case,
\begin{align*}
&\mathbcal{h}(\Phi(i_k\text{-block}))\\
&= \begin{cases}
				\mathbcal{h}(i_k\text{-block})+1, \quad &\text{if } \alpha_{0,1}=0~\text{and}~\displaystyle\sum_{k\geq 0}\left(2\alpha_{k,8} + \alpha_{k,4}\right)\geq 2,~\text{or if}~\alpha_{0,1}=1;\\
				\mathbcal{h}(i_k\text{-block}), \quad &\text{otherwise}.
	\end{cases}
\end{align*}
In this case also, there is no loss of 2-hooks under $\Phi$.\\
\noindent\textbf{Case 3: $\sum_{k\geq 0}\left(\alpha_{k,8} + \alpha_{k,4}\right)\neq0$ and $\sum_{k\geq 0}\left(\alpha_{k,5} + \alpha_{k,2}\right)\neq0$.}
In this case, if $k\neq 0$ or $i\neq 1$, then
\begin{align*}
	&\mathbcal{h}(\Phi(i_k\text{-block}))\\
	&= \begin{cases}
		\mathbcal{h}(i_k\text{-block})-1, \quad &\text{if both } \alpha_{k,4i} \text{ and } \alpha_{k,3i-1} \text{ are non-zero, where } i\in\{1,2\};\\
\mathbcal{h}(i_k\text{-block}), \quad &\text{otherwise}.
	\end{cases}
\end{align*}
If $k=0$ and $i=1$, then
\begin{align*}
	&\mathbcal{h}(\Phi(1_0\text{-block})) \\
	&=  \begin{cases}
		\mathbcal{h}(1_0\text{-block})-2, \quad &\text{if $1_0$-block} =\left(4,2^{\alpha_{0,2}}\right), \text{where } \alpha_{0,2}\neq 0, \displaystyle\sum_{k\geq 0}\left(\alpha_{k,8} + \alpha_{k,4}\right)=1;\\
		\mathbcal{h}(1_0\text{-block})-1, \quad &\text{either $1_0$-block}= (4^{\alpha_{0,4}},2^{\alpha_{0,2}}), \displaystyle\sum_{k\geq 0}\left(\alpha_{k,8} + \alpha_{k,4}\right)\geq2\\
        & \text{or $1_0$-block}= (4^{\alpha_{0,4}},2^{\alpha_{0,2}},1^{\alpha_{0,1}}),\alpha_{0,4}\neq0, \alpha_{0,2}\neq0,\alpha_{0,1}\geq2\\
        & \text{or $1_0$-block}= (2^{\alpha_{0,2}}), \alpha_{0,2}\neq0, \displaystyle\sum_{k\geq 0}\left(\alpha_{k,8} + \alpha_{k,4}\right)=1;\\
	\mathbcal{h}(1_0\text{-block})~\text{or}\quad &\text{otherwise.}\\ \mathbcal{h}(1_0\text{-block})+1,
	\end{cases}
\end{align*}
\par We see that under $\Phi$, only the partitions from \textbf{Case 3} may lose 2-hooks. We collect them in the following set: 
$$\mathbcal{C}_3(n) := \{\tau\in\mathbcal{B}_3(n) : \tau \text{ belongs to the partitions of \textbf{Case 3}}\}.$$
In our proof of Theorem \ref{Th1}, we compensate the loss of 2-hooks in the partitions from $\mathbcal{C}_3(n)$ by the partitions from $\mathbcal{B}_4(n)$ which are not in $\Phi(\mathcal{B}_3(n))$. We denote this set of partitions by
$$\mathbcal{C}_4(n) := \mathbcal{B}_4(n)\setminus \Phi(\mathcal{B}_3(n)).$$
\par We split $\mathbcal{C}_3(n)$ into the following two subsets and we handle them separately in Sections \ref{subsec2} and \ref{subsec1}, respectively:
\begin{itemize}
	\item $\mathbcal{C}_3^1(n):=\{\tau\in \mathbcal{C}_3(n) : \alpha_{0,1}\in \{0,1\}\}$,
		\item $\mathbcal{C}_3^2(n):=\mathbcal{C}_3(n) \setminus\mathbcal{C}_3^1(n)$.
\end{itemize}

Similarly, we split $\mathbcal{C}_4(n)$ going through a dissection of $\mathbcal{C}_4(n-2)$. We split $\mathbcal{C}_4(n-2)$ into the following three disjoint subsets.
\begin{itemize}
	\item $\mathbcal{Q}_1:=\{\tau\in\mathbcal{C}_4(n-2) : \alpha_{0,1}=0\}$,
	\item $\mathbcal{Q}_2:=\{\tau\in\mathbcal{C}_4(n-2) : \alpha_{0,1}\neq0 \text{ and } \sum_{k\geq 0}\alpha_{k,9}=0\}$,
	\item $\mathbcal{Q}_3:=\{\tau\in\mathbcal{C}_4(n-2) : \alpha_{0,1}\neq0 \text{ and } \sum_{k\geq 0}\alpha_{k,9}\neq 0\}$.
\end{itemize}
Let $(\tau,a)$ (resp. $(\tau,-a)$) denote the adjoining (resp. removing) of the part $a$ to (resp. from) the partition $\tau$. Now, we construct the following three subsets of $\mathbcal{C}_4(n)$.
\begin{itemize}
	\item $\mathbcal{Q}_1^\prime :=\{(\tau,2) : \tau\in\mathbcal{Q}_1\}$.
	\item  $\mathbcal{Q}_2^\prime :=
    \{(\tau,1,1) : \tau\in\mathbcal{Q}_2, \text{and } \alpha_{0,1}< \sum_{k\geq 0}(2\alpha_{k,6}+\alpha_{k,3})-2\}\\
		\null\hspace{11mm}\bigcup\{(\tau,2) : \tau\in\mathbcal{Q}_2, \text{and } \alpha_{0,1}=\sum_{k\geq 0}(2\alpha_{k,6}+\alpha_{k,3})-2\\
		\null\hspace{40mm} \text{ or } \alpha_{0,1}=\sum_{k\geq 0}(2\alpha_{k,6}+\alpha_{k,3})-1\}.$\\
 If $\alpha_{0,1}\geq \sum_{k\geq 0}(2\alpha_{k,6}+\alpha_{k,3})$, then $\tau\notin \mathbcal{C}_4(n-2)$, because in this case $\tau$ is an image under the map $\Phi$. Furthermore, we cannot adjoin two 1's to $\tau$ if $\alpha_{0,1}=\sum_{k\geq 0}(2\alpha_{k,6}+\alpha_{k,3})-2$ or $\alpha_{0,1}=\sum_{k\geq 0}(2\alpha_{k,6}+\alpha_{k,3})-1$, since in these cases $(\tau,1,1)$ cannot be in $\mathbcal{C}_4(n)$.
	\item $\mathbcal{Q}_3^\prime :=\{(\tau,1,1) : \tau\in\mathbcal{Q}_3\}$. This transformation is possible because the images under the map $\Phi$ contain no parts of the from $12k+9$ for any $k\geq 0$.
\end{itemize}
Clearly, $\mathbcal{Q}_1^\prime,\mathbcal{Q}_2^\prime$, and $\mathbcal{Q}_3^\prime$ are disjoint. Let $\mathbcal{Q}:=\mathbcal{Q}_1^\prime\cup \mathbcal{Q}_2^\prime\cup \mathbcal{Q}_3^\prime$.
\par 
Our proof of Theorem \ref{Th1} is by induction. In Sections \ref{subsec1} and \ref{subsec2}, we show that, under the map $\Phi$, the loss of 2-hooks of the partitions in $\mathbcal{C}_3^2(n)$ and $\mathbcal{C}_3^1(n)$ is compensated by the 2-hooks of the partitions in $\mathbcal{Q}$ and $\mathbcal{C}_4(n)\setminus\mathbcal{Q}$, respectively.
\subsection{Compensation of the loss of 2-hooks of partitions in $\mathbcal{C}_3^2(n)$}\label{subsec1}
If possible, suppose that the loss of 2-hooks in the partitions in $\mathbcal{C}_3(n-2)$ under the map $\Phi$ can be compensated by the partitions in $\mathbcal{C}_4(n-2)$. Then there must be a surjective map from $\mathbcal{C}_3(n-2)$ to $\mathbcal{C}_4(n-2)$, and the map may be many-valued. 
Now, we show that the loss of 2-hooks of the partitions in $\mathbcal{C}_3^2(n)$ is compensated by the 2-hooks of the partitions in $\mathbcal{Q}$.
\par 
If we remove two 1's from the partitions in $\mathbcal{C}_3^2(n)$ then the set $\mathbcal{C}_3^2(n)$ transforms into $\mathbcal{C}_3(n-2)$. Reversely, if we adjoin two 1's in the partitions in $\mathbcal{C}_3(n-2)$, then we get all the partitions in $\mathbcal{C}_3^2(n)$. 
\par 
When we adjoin two 1's to a partition $\tau$ in $\mathbcal{C}_3(n-2)$, then, only the $1_0$-block is changed, but it does not lose any 2-hook. The remaining blocks of $\tau$ in $\mathbcal{C}_3(n-2)$ are identical to those of the partition $(\tau,1,1)$ in $\mathbcal{C}_3^2(n)$. Since we have assumed that the loss of 2-hooks in the partitions in $\mathbcal{C}_3(n-2)$ under the map $\Phi$ is compensated by the partitions in $\mathbcal{C}_4(n-2)$, the loss of 2-hooks in the blocks other than $1_0$-block of the partitions in $\mathbcal{C}_3^2(n)$ is compensated by $\mathbcal{Q}$. The loss of 2-hooks in $1_0$-blocks in the partitions in $\mathbcal{C}_3^2(n)$ will also be compensated by the $1_0$-blocks in $\mathbcal{Q}$, if they do not lose any 2-hook under the construction of $\mathbcal{Q}_1^\prime$, $\mathbcal{Q}_2^\prime$, and $ \mathbcal{Q}_3^\prime$.
\par 
Clearly, while constructing $\mathbcal{Q}_1^\prime$ and $ \mathbcal{Q}_3^\prime$, partitions do not lose any 2-hook in their $1_0$-blocks. Only the partitions of the form $\left(\ldots,3^{\alpha_{0,3}},1^{\alpha_{0,1}}\right)$ in $\mathbcal{Q}_2$ may lose a 2-hook in their $1_0$-block under the construction of $\mathbcal{Q}_1^\prime$, $\mathbcal{Q}_2^\prime$, and $ \mathbcal{Q}_3^\prime$. More precisely, when a partition $\tau$ in  $\mathbcal{Q}_2$ is of the form
\begin{equation}\label{form1}
	 \left(\ldots,3^{\alpha_{0,3}},1^{\alpha_{0,1}}\right), \text{ with } \alpha_{0,1}=\sum_{k\geq 0}(2\alpha_{k,6}+\alpha_{k,3})-2 \text{ or } \alpha_{0,1}=\sum_{k\geq 0}(2\alpha_{k,6}+\alpha_{k,3})-1,
	\end{equation}
then the $1_0$-block of $(\tau,2)$ loses one 2-hook. But, the $1_0$-block of the partitions of the following three forms in $\mathbcal{C}_4(n-2)$ gains one 2-hook under the construction of $\mathbcal{Q}_1^\prime$, $\mathbcal{Q}_2^\prime$, and $ \mathbcal{Q}_3^\prime$.
\begin{itemize}
	\item $\textbf{F1:} \left(\ldots,2^1\right)$,
	\item $\textbf{F2:} \left(\ldots,(12k+9)^{\alpha_{k,9}},\ldots,1^1\right)$,
	\item $\textbf{F3:} \left(\ldots,2^1,1^{\alpha_{0,1}}\right)$ in  $\mathbcal{Q}_2$, with $\alpha_{0,1}=\sum_{k\geq 0}(2\alpha_{k,6}+\alpha_{k,3})-2$ or $\alpha_{0,1}=\sum_{k\geq 0}(2\alpha_{k,6}+\alpha_{k,3})-1$.
\end{itemize}
In Section \ref{sec:appendix}, we prove that for each partition of the form \eqref{form1}, there corresponds a partition of one of these three forms, for $n-2\neq 7$. Hence, the loss of 2-hooks of the partitions in $\mathbcal{C}_3^2(n)$ is compensated by the 2-hooks of the partitions in $\mathbcal{Q}$, for all $n\geq 10$.
\subsection{Compensation of the loss of 2-hooks of partitions in $\mathbcal{C}_3^1(n)$}\label{subsec2}
In this section, we deal with the remaining subsets of $\mathbcal{B}_3(n)$ and $\mathbcal{B}_4(n)$, respectively:
\begin{itemize}
	\item $\mathbcal{D}_3(n) := \mathbcal{C}_3^1(n)$,
	\item $\mathbcal{D}_4(n) := \mathbcal{C}_4(n)\setminus \mathbcal{Q}$.
\end{itemize}
We divide the partitions in $\mathbcal{D}_4(n)$ into the following six disjoint classes. 
\begin{align*}
	\textbf{CL1:} & \ \{\tau\in\mathbcal{B}_4(n) : \alpha_{0,1}=0 \text{ and } \alpha_{0,2}=0\} \cup \{(\ldots,1^{\alpha_{0,1}})\in\mathbcal{C}_4(n) : \alpha_{0,1}= 1 \text{ or } 2\},\\
	&\text{ where } \sum_{k\geq 0}\alpha_{k,9}\neq 0,\\
	&\text{For the remaining classes, we have} \sum_{k\geq 0}\alpha_{k,9}= 0.\\
		\textbf{CL2:} & \ \{\tau\in\mathbcal{B}_4(n) : \alpha_{0,1}=0, \alpha_{0,2}=0, \text{ and } \sum_{k\geq 0}(\alpha_{k,6}+\alpha_{k,3})\neq 0\},\\
		\textbf{CL3:} & \ \{(\ldots,2^0,1)\in\mathbcal{C}_4(n) : \sum_{k\geq 0}(2\alpha_{k,6}+\alpha_{k,3})>1\},\\
		\textbf{CL4:} & \ \{(\ldots,2^0,1^2)\in\mathbcal{C}_4(n) : \sum_{k\geq 0}(2\alpha_{k,6}+\alpha_{k,3})>2\},\\
		\textbf{CL5:} & \ \{(\ldots,2^{\alpha_{0,2}},1)\in\mathbcal{C}_4(n) : \alpha_{0,2}\neq 0 \text{ and } \sum_{k\geq 0}(2\alpha_{k,6}+\alpha_{k,3})>3\},\\
		\textbf{CL6:} & \ \{(\ldots,2^{\alpha_{0,2}},1^2)\in\mathbcal{C}_4(n) : \alpha_{0,2}\neq 0 \text{ and } \sum_{k\geq 0}(2\alpha_{k,6}+\alpha_{k,3})>4\}.
\end{align*}
Now, we define a map $\Psi$ from $\mathbcal{D}_3(n)$ to $\mathbcal{D}_4(n)$, which transforms the parts of the form $12k+8$ (for $k\geq 0$) and $12k+4$ (for $k\geq 1$) of the partitions in $\mathbcal{D}_3(n)$ as follows.
\begin{align*}
	\Psi\left((12k+8)^{\alpha_{k,8}}\right) & = \begin{cases}
		(12k+9), &\text{ if } \alpha_{k,8}=1\\
		\left(12k+9, (12k+6)^{\alpha_{k,8}-1}\right), &\text{ otherwise }
	\end{cases}, &\text{ for } k\geq 0;\\
	\Psi\left((12k+4)^{\alpha_{k,4}}\right) & = (12k+3)^{\alpha_{k,4}}, &\text{ for } k\geq 1.
\end{align*} 
The other parts of size greater than 8 are mapped to themselves under $\Psi$. This process yields an integer $\mathbcal{r}$ in total, which may be negative if $\alpha_{k,8}=1$ for at least one $k\geq 0$ and it is given by
$$\mathbcal{r} =\displaystyle\sum_{\substack{k\geq0\\ \alpha_{k,8}=1}}(-1)+\sum_{\substack{k\geq0\\ \alpha_{k,8}>1}}(2\alpha_{k,8}-3)+\sum_{k\geq1}\alpha_{k,4}.$$
At certain instances, if required, we map parts of size 8 or more to some different parts unlike the above mappings. The action of $\Psi$ on the parts of sizes 1, 2, 4, 5, and 7 of $\tau\in\mathcal{D}_3(n)$ is yet to be defined. For $\tau\in\mathcal{D}_3(n)$, we consider the block $\left(7^{\alpha_{0,7}},5^{\alpha_{0,5}},4^{\alpha_{0,4}},2^{\alpha_{0,2}},1^{\alpha_{0,1}}\right)$ and we call it $\Xi$-block of $\tau$, and we call the remaining whole block of $\tau$ the $\Xi^\complement$-block. Let $\Xi:=7\times\alpha_{0,7}+5\times\alpha_{0,5}+4\times\alpha_{0,4}+2\times\alpha_{0,2}+1\times\alpha_{0,1}$. Moreover, we define the $\overrightarrow{\Xi}$-block in $\mathbcal{D}_4(n)$ as $\left(7^{\alpha_{0,7}},6^{\alpha_{0,6}},5^{\alpha_{0,5}},3^{\alpha_{0,3}},2^{\alpha_{0,2}},1^{\alpha_{0,1}}\right)$, and we define $\overrightarrow{\Xi}$ as $7\times\alpha_{0,7}+6\times\alpha_{0,6}+5\times\alpha_{0,5}+3\times\alpha_{0,3}+2\times\alpha_{0,2}+1\times\alpha_{0,1}$. Note that when $\alpha_{0,8}>1$ in the $\Xi^{\complement}$-block then the occurrence of part 6 in the image partitions is not considered in the $\overrightarrow{\Xi}$-block.
\par 
Under the map $\Phi$, $3_k$-block of a partition never loses a 2-hook ($k\geq0$); $2_k$-block may lose one 2-hook in the case $\alpha_{k,8}\neq0$ ($k\geq0$); $1_k$-block may lose one 2-hook in the case $\alpha_{k,4}\neq0$ ($k\geq1$). However, under the map $\Psi$, when $\alpha_{k,8}\neq 0$ with $k\geq0$, then $\Psi\left((12k+8)^{\alpha_{k,8}}\right)$ yields at least one 2-hook in the $3_k$-block. For $k\geq 1$, the transformation $\Psi\left((12k+4)^{\alpha_{k,4}}\right)$ yields at least one 2-hook in the $1_k$-block when $\alpha_{k,4}\neq 0$, since $12k$ is not a part. Therefore, to compensate the loss of 2-hooks of partitions under the map $\Phi$ through the map $\Psi$, our concern is the $1_0$-block.
\par
We know that if the $1_0$-block is of the form $(4,2^{\alpha_{0,2}})$, where $\alpha_{0,2}\neq 0$, then it may lose two 2-hooks under the map $\Phi$, otherwise it loses at most one 2-hook. Therefore, when the $\Xi$-block is of the form $(4,2^{\alpha_{0,2}})$, where $\alpha_{0,2}\neq 0$, and $\mathbcal{r}=2s-1$, where $s$ is a nonnegative integer, then the block can be transformed as follows to obtain two 2-hooks.
\begin{align*}
	\Psi((4,2^{\alpha_{0,2}}))= \begin{cases}
		\left(3^{\frac{2(\alpha_{0,2}+s+2)}{3}-1},1^2\right), & \text{ if } \alpha_{0,2}+s\equiv1\pmod3;\\
\left(3^{\frac{2(\alpha_{0,2}+s+1)}{3}},1\right), & \text{ if } \alpha_{0,2}+s\equiv2\pmod3;\\
\left(3^{\frac{2(\alpha_{0,2}+s)}{3}+1}\right), & \text{ if } \alpha_{0,2}+s\equiv0\pmod3.
	\end{cases}
\end{align*}
Here, $\overrightarrow{\Xi}$ is odd and greater than or equal to 5. Similarly, when $\mathbcal{r}=2s$, where $s$ is a nonnegative integer, we can transform as follows.
\begin{align*}
	\Psi((4,2^{\alpha_{0,2}}))= \begin{cases}
		\left(3^{\frac{2(\alpha_{0,2}+s+2)}{3}}\right), & \text{ if } \alpha_{0,2}+s\equiv1\pmod3;\\
		\left(5, 3^{\frac{2(\alpha_{0,2}+s+1)}{3}-1}\right), & \text{ if } \alpha_{0,2}+s\equiv2\pmod3;\\
		\left(3^{\frac{2(\alpha_{0,2}+s)}{3}+1},1\right), & \text{ if } \alpha_{0,2}+s\equiv0\pmod3.
	\end{cases}
\end{align*}
Here, $\overrightarrow{\Xi}$ is even and greater than or equal to 6. For example, when the $\Xi^\complement$-block for the partitions of $n=20$ is $(10)$, then we get a subset $\{(10, 5, 4, 1), (10, 4^2, 2), (10, 4, 2^3)\}$ of $\mathbcal{D}_3(20)$. When $\mathbcal{r}=0$, to apply $\Psi$, the corresponding subset of $\mathbcal{D}_4(20)$ is $\{(10, 7, 3), (10, 6, 3, 1), (10, 3^3, 1)\}$. But the partition $(10, 4, 2^3)$ loses two 2-hooks under $\Phi$. Therefore, we take $\Psi((10,4,2^3))= (10,3^3,1)$. Again, when $\mathbcal{r}=1$, $n=22$, and $\Xi^\complement$-block is $(16)$, then for
$\{(16, 5, 1), (16, 4, 2), (16, 2^3)\} \subset \mathbcal{D}_3(22)$, we have the corresponding subset $\{(15, 7), (15, 6, 1), (15, 3^2, 1)\}$ in $ \mathbcal{D}_4(22)$. Here, we take $\Psi((16, 4, 2))=(15, 3^2, 1)$.
\par
Furthermore, when the $\overrightarrow{\Xi}$-block of a partition in $\mathbcal{D}_4(n)$ is $(3,2,1)$, then the $1_0$-block of the partition has no 2-hook. When necessary, we map the partition in $\mathbcal{D}_4(n)$ from a partition in $\mathbcal{D}_3(n)$ whose $1_0$-block has no loss of 2-hooks under $\Phi$. When $\Xi=1,2,3,4,5,6,$ or $7$, then this $\overrightarrow{\Xi}$-block may arise, and for each of these values we have such partitions in $\mathbcal{D}_3(n)$. For instance, if $\Xi=5$, then the $1_0$-block $(2^2,1)$ in $\mathbcal{D}_3(n)$ has no loss of 2-hooks under $\Phi$. Similarly, when $\Xi=1,2,3$, or $4$, then the $\overrightarrow{\Xi}$-block can be $(2,1)$, and when $\Xi=1$, or $2$, then the $\overrightarrow{\Xi}$-block can be $(1)$. In both the cases, the $1_0$-block in the pre-images lose no 2-hooks. Moreover, if the smallest part in a $\overrightarrow{\Xi}$-block is greater than 1, then it has at least one 2-hook.
\par 
Now, we are in a position to define $\Psi(\Xi\text{-block})$ for the remaining types of $\Xi$-blocks. In the following cases, for a fix $\Xi^\complement$-block for $n$, we find the generating functions for the number of $\Xi$-blocks and the number of $\overrightarrow{\Xi}$-blocks corresponding to $\Psi(\Xi^\complement\text{-block})$. We denote them by $\mathfrak{X}_3(q)$ and $\mathfrak{X}_4(q)$, respectively. Clearly, the nonnegativity of the coefficients of $q^n$ in $\mathfrak{X}_4(q)-\mathfrak{X}_3(q)$ implies that the number of required partitions in $\mathbcal{D}_4(n)$ is greater or equal to the number of partitions in $\mathbcal{D}_3(n)$ in the following cases. Also, if a partition $\tau\in\mathbcal{D}_3(n)$ has a loss of 2-hook under $\Phi$, then it will be compensated by the 2-hook of $\Psi(\tau)$.\\
\noindent \textbf{Case A: $\sum_{k\geq 0}\alpha_{k,8}\neq 0$.} We divide this case into the following three subcases.\\
\noindent \textit{\textbf{Case A1:} $\mathbcal{r}=-1$.}\\
\noindent \textit{\textbf{Case A1a:} $\sum_{k\geq 1}\alpha_{k,4}=0$.} We map the partitions from $\mathbcal{D}_3(n)$ falling in this case to the partitions of class $\textbf{CL1}$ of $\mathbcal{D}_4(n)$ and we show that there are more such partitions in $\mathbcal{D}_4(n)$ than in $\mathbcal{D}_3(n)$. Here, 6 is not a part in the corresponding $\overrightarrow{\Xi}$-blocks. Moreover, the image partitions contain a part of the form $12k+9$. In this case, we have $\mathfrak{X}_4(q)$ corresponding to $\Xi+\mathbcal{r}$ as follows,
\begin{equation*}
	\mathfrak{X}_4(q)=q\left(\frac{1}{(1-q^3)(1-q^5)(1-q^7)}+\frac{q^1+q^{1+1}}{(1-q^2)(1-q^3)(1-q^5)(1-q^7)}\right).
\end{equation*}
Also,
\begin{equation*}
	\mathfrak{X}_3(q)=\frac{1+q^1}{(1-q^2)(1-q^4)(1-q^5)(1-q^7)}.
\end{equation*}
Therefore,
\begin{equation*}
	\mathfrak{X}_4(q)-\mathfrak{X}_3(q)=\frac{-1+q^3+q^4}{(1-q^3)(1-q^4)(1-q^5)(1-q^7)}.
\end{equation*}
Next, we check the nonnegativity of the coefficients of $\mathfrak{X}_4(q)-\mathfrak{X}_3(q)$. Let $P_a(n)$ be the set of all partitions of $n$ with parts $3,~4,~5$, or $7$ only, and let $p_a(n)$ be its cardinality. Then, $\mathfrak{X}_4(q)-\mathfrak{X}_3(q)$ is the generating function for $-p_a(n)+p_a(n-3)+p_a(n-4)$. In the later cases as well, we use the same notation for the partition statistics corresponding to $\mathfrak{X}_4(q)-\mathfrak{X}_3(q)$. 
\par 
Now, we adjoin 3 as a part to the partitions in $P_a(n-3)$ to obtain all the partitions in $P_a(n)$ that have at least one 3 as a part. We remove these partitions from $P_a(n)$.
\par 
Next, we transform each remaining partition in $P_a(n)$ into a partition in $P_a(n-4)$ as follows.
\begin{align*}
	(\ldots,4)\to & \ 	(\ldots), \text{ removed the smallest part}\\
	(\ldots,5,5)\to & \	(\ldots,3^2)\\
	(\ldots,7,5)\to & \	(\ldots,5,3^1)\\
	(\ldots,7,7)\to & \	(\ldots,7,3^1).
\end{align*}
Each of these types of partitions in $P_a(n-4)$ is distinct. So, we get
\[-p_a(n)+p_a(n-3)+p_a(n-4)\geq 0, \text{ for all } n\geq 7.\]
Checking for the lower values of $n$, we obtain that the coefficients of $\mathfrak{X}_4(q)-\mathfrak{X}_3(q)$ are nonnegative for all $n\geq 1$, except when $n=5$. Therefore, when $\Xi=5$, some partitions of this category cannot be mapped by $\Psi$ to distinct partitions in $\mathbcal{D}_4(n)$.
\par 
When $\Xi=5$ for a fix $\Xi^\complement$-block, the $\Xi$-blocks are $(5)$, $(4,1)$, and $(2,2,1)$, and the $\overrightarrow{\Xi}$-blocks are $(3,1)$ and $(2,1,1)$. We see that for $k> 1$, if the $2_k$-block of a partition does not lose a 2-hook under $\Phi$, then under $\Psi$ the $2_k$-block and $3_k$-block have at least one 2-hook in total. Again, if the $2_k$-block loses a 2-hook under $\Phi$, then after compensating this 2-hook by applying $\Psi$ we still have at least one 2-hook in the $2_k$-block and $3_k$-block in total, for all $k> 1$. Since, the $1_0$-blocks lose no 2-hooks under the map $\Phi$, in the similar manner, we get 2-hooks under $\Psi$ for $k=1$ as well. Therefore, the total loss of 2-hooks in the three partitions in $\mathbcal{D}_3(n)$ is compensated by the total 2-hooks in the two partitions in  $\mathbcal{D}_4(n)$. 
\par 
When $\Xi=0$, a partition $\tau$ may lose 2-hooks under $\Phi$ if it has at least one part of the form $12k+5$ for $k\geq 1$. Then, we transform the smallest of such parts of $\tau$ as follows.
\[\left((12k+5)^{\alpha_{k,5}}\right) \to \left((12k+5)^{\alpha_{k,5}-1},12(k-1)+7, 9\right).\]
For example, $(20,17^2)\to (21,17,9,7)$, and $(32,20,17,8^2)\to (33,21,9^2,7,6)$. When $8$ is a part of $\tau$, then $9^2$ is a part of its image. Therefore, this image is different from the other partitions in $\mathbcal{D}_4(n)$ that have already been mapped. Moreover, when $8$ is not a part of $\tau$, then the image can not have a pre-image under the general transformation method of the map $\Psi$.
\par 
Thus, the loss of 2-hooks under $\Phi$ in the partitions $$\left(\Xi^\complement\text{-block},~\Xi\text{-block}\right)$$ is compensated by the 2-hooks in the partitions $$\left(\Psi\left(\Xi^\complement\text{-block}\right),~\overrightarrow{\Xi}\text{-block}\right).$$
\noindent \textit{\textbf{Case A1b:} $\sum_{k\geq 1}\alpha_{k,4}\neq 0$.}
The arguments in this case are similar to those in \textit{\textbf{Case A1a}}, with the exception of one more transformation each for $\Xi=0$, and $5$.
\par 
When $\Xi=5$, similar to \textit{\textbf{Case A1a}}, the total loss of 2-hooks in the three partitions in $\mathbcal{D}_3(n)$ is compensated by the total 2-hooks in the two partitions in  $\mathbcal{D}_4(n)$, since $\sum_{k\geq 0}\alpha_{k,8}>\sum_{k\geq 1}\alpha_{k,4}$.
\par 
When $\Xi=0$, and the partition $\tau$ has at least one part of the form $12k+5$ for $k\geq 1$, then we transform it similarly as done in \textit{\textbf{Case A1a}}. If there is no part of the form $12k+5$, then there must be at least one part of the form $12k+2$. We then transform the smallest of these parts of the form $12k+2$ as follows.
\begin{align*}
	\left((12k+2)^{\alpha_{k,2}}\right) \to\begin{cases}
		(9,3,1), & \text{if the part is } 14;\\
		\left((12k+2)^{\alpha_{k,2}-1},12(k-2)+11, 9,3,1^2\right), & \text{otherwise}.
	\end{cases}
\end{align*}
For example,
\begin{align*}
	(44,32,28,20^2,14,8) & \to (45,33,27,21,18,9^2,3,1),\\
	(56,44,32,28,20^2,14) & \to (57,45,33,27,21,18,9,3,1),\\
	(44,32,28^2,26^2,20) & \to (45,33,27^2,26,21,11,9,3,1^2).
\end{align*}

\noindent \textit{\textbf{Case A2:} $\mathbcal{r}> -1$.}  Here, 
\begin{equation*}
	\mathfrak{X}_4(q)=q^{-\mathbcal{r}}\left(\frac{1}{(1-q^3)(1-q^5)(1-q^7)}+\frac{q^1+q^{1+1}}{(1-q^2)(1-q^3)(1-q^5)(1-q^7)}\right).
\end{equation*}
Therefore, we have that for each $\mathbcal{r}$, the coefficients of $\mathfrak{X}_4(q)-\mathfrak{X}_3(q)$ are nonnegative for all $n\geq 0$.
\\
\noindent\textit{\textbf{Case A3:} $\mathbcal{r}\leq -2$.}
If the largest part of the form $(12k+8)^1$ is $12\ell+8$, then $\mathbcal{r}\geq -(\ell+1)$. In this case, we transform the part $12\ell+8$ in a different way as follows.
Let
\[12\ell+8+\mathbcal{r}+1=9a+b,\]
for some nonnegative integers $a$ and $b$, where $0\leq b <9$. Then, $(12\ell+8)\to (9^a)$,
and the $\overrightarrow{\Xi}$-block must be derived from $\Xi+b$. In this process we have, $\alpha_{0,9}$ is at least 3. For example,
\[(47,44,32,20,16^2,14,8,4,2)\to (47,33,21,15^2,14,9^5,(\overrightarrow{\Xi}\text{-block})),\]
where we must derive $\overrightarrow{\Xi}$-block from 13, since for this partition $12\ell+8=44$, $\mathbcal{r}=-2$, $a=4$, and $b=7$.
\par 
In this case, when the $\Xi$-block is of the form $(4,2^{\alpha_{0,2}})$, where $\alpha_{0,2}\neq 0$, then we transform the block as mentioned above by taking $\mathbcal{r}=b$, and $s=\big\lceil\frac{b}{2}\big\rceil$. For instance, for the above example $(4,2)\to (3^4,1)$.
\begin{remark}
The image partitions in all three cases \textbf{Case A1}, \textbf{Case A2}, and \textbf{Case A3} are from class \textbf{CL1} and are distinct.  
\end{remark}
\noindent \textbf{Case B: $\sum_{k\geq 0}\alpha_{k,8}=0$.}
Here, 6 may be a part in the corresponding $\overrightarrow{\Xi}$-blocks. We have the following two subcases.\\
\noindent\textit{\textbf{Case B1:} $\mathbcal{r}=0$.}
That is, $\sum_{k\geq 1}\alpha_{k,4}=0$. Therefore, 4 must be a part in $\Xi$-block, and 3 or 6 must be a part in $\overrightarrow{\Xi}$-block. In view of classes \textbf{CL2}, \textbf{CL3}, \textbf{CL4}, \textbf{CL5}, and \textbf{CL6}, we have
\begin{multline*}
	\mathfrak{X}_4(q)= \frac{1}{(1-q^3)(1-q^5)(1-q^6)(1-q^7)}-\frac{1}{(1-q^5)(1-q^7)}\\
	+\frac{q^1+q^{1+1}}{(1-q^2)(1-q^3)(1-q^5)(1-q^6)(1-q^7)} - \frac{q^1+q^{1+3}+q^{1+2}(q^{3+3}+q^{3+3+3}+q^{3+6}+q^{6})}{(1-q^2)(1-q^5)(1-q^7)}\\
		-\frac{q^{1+1}+q^{1+1}(q^{3}+q^{3+3}+q^{6})+q^{1+1+2}(q^{3+3+3}+q^{3+3+3+3}+q^{3+6}+q^{3+3+6}+q^{6+6})}{(1-q^2)(1-q^5)(1-q^7)}.
\end{multline*}
In this case, if 2 is not a part in the partitions in $\mathbcal{D}_3(n)$, then there is no loss of 2-hooks. Therefore, we ignore those partitions. We derive the generating function for the number of partitions with 2 as a part. Hence,
\begin{equation*}
	\mathfrak{X}_3(q)= \frac{q^{2+4}+q^{1+2+4}}{(1-q^2)(1-q^4)(1-q^5)(1-q^7)}.
\end{equation*}
Thus, $\mathfrak{X}_4(q)-\mathfrak{X}_3(q)$ becomes a fraction whose numerator and denominator are $q^3(1 - q^5 - q^6 - q^7 + q^9 + q^{10} + 2q^{11} + q^{12} + 2q^{14} - q^{15} + 3q^{16} + 3q^{18} + q^{19} - 2q^{20} + q^{21} - 3q^{22} - 3q^{24})$ and $(1-q^3)(1-q^4)(1-q^5)(1-q^6)(1-q^7)$, respectively.
Let $p_b(n)$ denote the number of partitions of $n$ with parts $3,4,5,6$, or $7$ only. Then, $\dfrac{1}{q^3}\left(\mathfrak{X}_4(q)-\mathfrak{X}_3(q)\right)$ is the generating function for $p_b(n)-p_b(n-5)-p_b(n-6)-p_b(n-7)+p_b(n-9)+p_b(n-10)+2p_b(n-11)+p_b(n-12)+2p_b(n-14)-p_b(n-15)+3p_b(n-16)+3p_b(n-18)+p_b(n-19)-2p_b(n-20)+p_b(n-21)-3p_b(n-22)-3p_b(n-24)$, which is greater than or equal to $p_b(n)-p_b(n-5)-p_b(n-6)-p_b(n-7)+p_b(n-9)+p_b(n-10)+p_b(n-11)+p_b(n-12)+2p_b(n-14)$. By a similar method to \textit{\textbf{Case A1}}, it can be easily proved that $-p_b(n-6)+p_b(n-9)+p_b(n-11)$ and $-p_b(n-7)+p_b(n-10)+p_b(n-12)$ are nonnegative for all $n\geq 20$, and $p_b(n)-p_b(n-5)$ is nonnegative for all $n\geq 0$, and finally, the coefficients of $\mathfrak{X}_4(q)-\mathfrak{X}_3(q)$ are nonnegative for all $n\geq 0$.\\
\noindent\textit{\textbf{Case B2:} $\mathbcal{r}\neq 0$.}
That is, $\sum_{k\geq 1}\alpha_{k,4}\neq 0$. We have,
\begin{equation*}
	\mathfrak{X}_3(q) = \frac{1+q}{(1-q^2)(1-q^4)(1-q^5)(1-q^7)}.
\end{equation*}
Here, 3 or 6 may not be a part in the $\overrightarrow{\Xi}$-block. Therefore, considering the following four subcases we transform the partitions. Similar to \textit{\textbf{Case B1}}, in these cases as well, the image partitions are from classes \textbf{CL2}, \textbf{CL3}, \textbf{CL4}, \textbf{CL5}, and \textbf{CL6}.\\
\noindent\textit{\textbf{Case B2a:} $\mathbcal{r}=1$.}
We have,
\begin{multline*}
	\mathfrak{X}_4(q)= \frac{1}{q}\bigg(\frac{1}{(1-q^3)(1-q^5)(1-q^6)(1-q^7)}
	+\frac{q^1+q^{1+1}}{(1-q^2)(1-q^3)(1-q^5)(1-q^6)(1-q^7)}\\
	-\frac{q^1+q^{1+2}(q^3+q^{3+3}+q^6)}{(1-q^2)(1-q^5)(1-q^7)}
	-\frac{q^{1+1}+q^{1+1+3}+q^{1+1+2}(q^{3+3}+q^{3+3+3}+q^{3+6}+q^6)}{(1-q^2)(1-q^5)(1-q^7)}\bigg).
\end{multline*}
Here, we find $\mathfrak{X}_4(q)-\mathfrak{X}_3(q)$, and observe that for $\Xi=0,1,4,5,7,8,9,10$, and $12$, some partitions of this category cannot be mapped by $\Psi$ to distinct partitions in $\mathbcal{D}_4(n)$. For $\Xi=0$, and $1$, we apply the following transformation,
\[\left(12k+4, (12k+2)^{\alpha_{k,2}}\right) \to \left((12k+3)^2, (12k+2)^{\alpha_{k,2}-1}\right).\]
For example, $(16,14^5)\to(15^2,14^4)$, and $(16,14^5,1)\to(15^2,14^4,1)$. Furthermore, for each $\Xi=4,5,7,8,9,10$, and $12$, we see that the total loss of 2-hooks in the partitions in $\mathbcal{D}_3(n)$ is smaller than the total 2-hooks in the corresponding partitions in $\mathbcal{D}_4(n)$.\\
\noindent\textit{\textbf{Case B2b:} $\mathbcal{r}=2$.}
We have,
\begin{multline*}
	\mathfrak{X}_4(q)= \frac{1}{q^2}\bigg(\frac{1}{(1-q^3)(1-q^5)(1-q^6)(1-q^7)}
	+\frac{q^1+q^{1+1}}{(1-q^2)(1-q^3)(1-q^5)(1-q^6)(1-q^7)}\\
	-\frac{q^{1+2+3}}{(1-q^2)(1-q^5)(1-q^7)}
	-\frac{q^{1+1}+q^{1+1+2}(q^3+q^{3+3}+q^6)}{(1-q^2)(1-q^5)(1-q^7)}\bigg).
\end{multline*}
We find $\mathfrak{X}_4(q)-\mathfrak{X}_3(q)$, and observe that the coefficient is negative for $\Xi=0$ only. For $\Xi=0$, we transform the part $12k+4$ if its rightmost $12k+2$ is also a part. Here, ${\alpha_{k,4}}$ is at most 2. The transformation is
\[\left((12k+4)^{\alpha_{k,4}}, (12k+2)^{\alpha_{k,2}}\right) \to \left((12k+3)^{\alpha_{k,4}+1}, (12k+2)^{\alpha_{k,2}-1}\right).\]
For example,
$(28,26,16)\to (27^2,15,1)$,
$(28,26,16,14)\to (27^2,15^2)$,
$(16^2,14)\to (15^3,1)$,
and 
$(28,16,14^9)\to (27,15^2,14^8,1)$.\\
\noindent\textit{\textbf{Case B2c:} $\mathbcal{r}=3$.}
We have,
\begin{multline*}
	\mathfrak{X}_4(q)= \frac{1}{q^3}\bigg(\frac{1}{(1-q^3)(1-q^5)(1-q^6)(1-q^7)}\\
	+\frac{q^1+q^{1+1}}{(1-q^2)(1-q^3)(1-q^5)(1-q^6)(1-q^7)}
	-\frac{q^{1+1+2}+q^{1+1+2+3}}{(1-q^2)(1-q^5)(1-q^7)}\bigg).
\end{multline*}
Here, the coefficients of $\mathfrak{X}_4(q)-\mathfrak{X}_3(q)$ are nonnegative for all $n\geq 0$.\\
\noindent\textit{\textbf{Case B2d:} $\mathbcal{r}\geq 4$.}
We have,
\begin{equation*}
	\mathfrak{X}_4(q)= \frac{1}{q^m}\bigg(\frac{1+q}{(1-q^2)(1-q^3)(1-q^5)(1-q^6)(1-q^7)}\bigg), \text{ where } m\geq 4.
\end{equation*}
In this case, the coefficient of $q^n$ in $\mathfrak{X}_4(q)-\mathfrak{X}_3(q)$ is given by
$p_b(n+m)+p_b(n+m-1)+p_b(n+m-2)+p_b(n+m-3)
-p_b(n)-p_b(n-1)-p_b(n-2)+p_b(n-6)+p_b(n-7)+p_b(n-8),$
which is positive. Therefore, for all $m\geq 4$, the coefficients of $\mathfrak{X}_4(q)-\mathfrak{X}_3(q)$ are positive for all $n\geq 0$. Table \ref{table_main} demonstrates the above mapping of partitions from $\mathbcal{D_3}(n)$ to partitions from $\mathbcal{D}_4(n)$, when $n=22$ and $\Xi^{\complement}$-block is nonzero.
\par 
Assuming for $\mathbcal{C}_3(n-2)$, we finally deduce for $\mathbcal{C}_3(n)$. Since $n$ is arbitrary, and it is true for the initial values of $n$, we conclude by induction that for all $n\geq0$, $b_{4,2}(n)\geq b_{3,2}(n)$. This completes the proof of Theorem \ref{Th1}.
\begin{center}
	\begin{table}
		\begin{tabular}[h]{|m{1.33cm}|m{3.7cm}|m{1,95cm}|m{4cm}|}
			\hline
			$\Xi^\complement$-block &
			$(\Xi^\complement\text{-block},\Xi\text{-block})$ &
			$\Psi(\Xi^\complement\text{-block})$ &
			$(\Psi(\Xi^\complement\text{-block}),\overrightarrow{\Xi}\text{-block})$ \\
			\hline
			$(20)$ &
			$(20,2)$ &
			$(21)$ & 
			$(21,1)$
			\\
			\hline
			$(16)$ &
			$(16,5,1)$, $(16,4,2)$, $(16,2^3)$ &
			$(15)$ &
			$(15,7)$, $(15,6,1)$, $(15,3^2,1)$\\
			\hline
			$(14,8)$ &
			$(14,8)$ &
			&\\
			\hline
			$(14)$ &
			$(14,4^2)$, $(14,4,2^2)$ &
			$(14)$ &
			$(14,5,3)$ \\
			\hline
			$(13)$ &
			$(13,5,4)$, $(13,4,2^2,1)$ &
			$(13)$ &
			$(13,6,3)$, $(13,3^3)$ \\
			\hline
			$(11,8)$ &
			$(11,8,2,1)$ &
			$(11,9)$ &
			$(11,9,1^2)$ 
			\\
			\hline
			$(11)$ &
			$(11,5,4,2)$, $(11,4^2,2,1)$, $(11,4,2^3,1)$ &
			$(11)$ &
			$(11,6,5)$, $(11,6,3,1^2)$, $(11,5,3^2)$, $(11,3^3,1^2)$ 
			\\
			\hline
			$(10,8)$ &
			$(10,8,2^2)$ &
			$(10,9)$ &
			$(10,9,3)$, $(10,9,2,1)$
			\\
			\hline
			$(10)$ &
			$(10,5,4,2,1)$, $(10,4^2,2^2)$, $(10,4,2^4)$&
			$(10)$ &
			$(10,6^2)$, $(10,6,5,1)$, $(10,6,3^2)$, $(10,5,3^2,1)$, $(10,3^4)$\\
			\hline
			$(8^2)$ &
			$(8^2,5,1)$, $(8^2,4,2)$, $(8^2,2^3)$&
			$(9,6)$&
			$(9,7,6)$, $(9,6,5,1^2)$, $(9,6,3^2,1)$, $(9,6,3,2,1^2)$, $(9,6,2^3,1)$ \\
			\hline
			$(8)$ &
			$(8,7,5,2)$, $(8,7,4,2,1)$, $(8,7,2^3,1)$, $(8,5^2,4)$, $(8,5^2,2^2)$, $(8,5,4^2,1)$, $(8,5,4,2^2,1)$, $(8,5,2^4,1)$, $(8,4^3,2)$, $(8,4^2,2^3)$, $(8,4,2^5)$, $(8,2^7)$ &
			$(9)$ &
			$(9,7,5,1)$, $(9,7,3^2)$, $(9,7,3,2,1)$, $(9,7,2^2,1^2)$, $(9,5^2,3)$, $(9,5^2,2,1)$, $(9,5,3^2,1^2)$, $(9,5,3,2^2,1)$, $(9,5,2^3,1^2)$, $(9,3^4,1)$, $(9,3^3,2,1^2)$, $(9,3^2,2^3,1)$, $(9,3,2^4,1^2)$
			\\
			\hline
			&&&
			$(19,3)$, $(18,3,1)$, $(13,9)$, $(9^2,2,1^2)$,  $(9^2,3,1)$, $(9,6^2,1)$ \\
			\hline
		\end{tabular}
		\caption{Mapping from $\mathcal{D}_3(22)$ to $\mathcal{D}_4(22)$, where the largest part of a partition in $\mathcal{D}_3(22)$ is greater than 7}
		\label{table_main}
	\end{table}
\end{center}
\section{Appendix}\label{sec:appendix}
Let $\tau$ be a partition in $\mathbcal{Q}_2$ of the form
$\left(\ldots,3^{\alpha_{0,3}},1^{\alpha_{0,1}}\right),$ with $\alpha_{0,1}=\sum_{k\geq 0}(2\alpha_{k,6}+\alpha_{k,3})-2 \text{ or } \alpha_{0,1}=\sum_{k\geq 0}(2\alpha_{k,6}+\alpha_{k,3})-1$. In this section, we show that for each such $\tau$, there exists a partition from one of the forms \textbf{F1}, \textbf{F2}, or \textbf{F3}. Depending on the multiplicity of part 1 in $\tau$, we consider the following four cases.

\noindent\textbf{Case (i). $\alpha_{0,1}=1$.} Then $\sum_{k\geq 0}(2\alpha_{k,6}+\alpha_{k,3})$ must be 3 or 2.

\noindent\textbf{Case (i-a). $\sum_{k\geq 0}(2\alpha_{k,6}+\alpha_{k,3})=3$.}
In this case, we correspond the partitions of the form \eqref{form1} to the partitions of the form \textbf{F1}, \textbf{F2}, or \textbf{F3} as follows:
\begin{align*}
	\left(\ldots,3^{\alpha_{0,3}},1\right) \to\begin{cases}
		\left(\ldots,12k+9,\ldots,1\right),
		&\text{if } 12k+6 \text{ is a part for a certain } k\geq 0;\\
		\left(\ldots,12k^\prime+5,\ldots,2\right),
		&\text{if } \sum_{k\geq 0}\alpha_{k,6}=0, \alpha_{0,3}=1, \text{ and $k^\prime$ is the  }\\
		&\text{smallest $k\geq1$ for which $\alpha_{k,3}\neq0$};\\ 
		\left(\ldots,5,2\right),
		& \text{if } \alpha_{0,3}=2;\\
		\left(\ldots,5,3,2\right),
		& \text{if } \alpha_{0,3}=3.
	\end{cases}
\end{align*}
Note that for the second, third, and forth partitions on the right-hand side of above correspondence, we have $\sum_{k\geq 0}(2\alpha_{k,6}+\alpha_{k,3})=1$.\\
\noindent\textbf{Case (i-b). $\sum_{k\geq 0}(2\alpha_{k,6}+\alpha_{k,3})=2$.}
Here, $\alpha_{0,3}=2$, or apart from one occurrence of the part 3 we have another part of the form $12k_1+3$ for some $k_1\geq 1$. The remaining parts can be of the form $12k+5$, $12k+7$, $12k+10$, and $12k+11$, where $k\geq 0$, and of the form $12k+1$, and $12k+2$, where $k\geq 1$. We consider the smallest part greater than 3 and not of the form $12k_1+3$. For simplicity, we denote it by $k_s$. Now, we have the following correspondence:
	\begin{align*}
		\left(\ldots,3^{\alpha_{0,3}},1\right) \to\begin{cases}
			&\left(\ldots,(12k_2+11)^{\alpha_{k_2,11}-1},12k_2+9,\ldots,2^{6\delta k_1+4}\ldots,1\right),\\
			& \qquad\qquad\qquad \text{if } k_s=12k_2+11;\\
			&\left(\ldots,(12k_2+10)^{\alpha_{k_2,10}-1},12k_2+9,\ldots,3,2^{6\delta k_1+2}\ldots,1\right),\\
			& \qquad\qquad\qquad \text{if } k_s=12k_2+10;\\
			&\left(\ldots,12k_2+9,(12k_2+7)^{\alpha_{k_2,7}-1},\ldots,2^{6\delta k_1+2}\ldots,1\right),\\
			& \qquad\qquad\qquad\text{if } k_s=12k_2+7;\\
			&\left(\ldots,12k_2+9,\ldots,(12k_2+5)^{\alpha_{k_2,5}-1},\ldots,2^{6\delta k_1+1}\ldots,1\right),\\
			& \qquad\qquad\qquad \text{if } k_s=12k_2+5;\\
			&\left(\ldots,(12k_2+2)^{\alpha_{k_2,2}-1},12(k_2-1)+9,\ldots,3,2^{6\delta k_1+4}\ldots,1\right),\\
			&  \qquad\qquad\qquad\text{if } k_s=12k_2+2;\\
			&\left(\ldots,(12k_2+1)^{\alpha_{k_2,1}-1},12(k_2-1)+9,\ldots,2^{6\delta k_1+5}\ldots,1\right),\\
			&  \qquad\qquad\qquad\text{if } k_s=12k_2+1;\\
			&\left(12(k_1-1)+9,3,2^3,1\right),\\
			&\qquad\qquad\qquad\text{if 1, 3, and $12k_1+3$ are the only parts};
		\end{cases}
	\end{align*}
	where, \begin{align*}
		\delta=\begin{cases}
			1 & \text{if } \alpha_{0,3}=1,\\
			0 & \text{if } \alpha_{0,3}=2.
		\end{cases}
\end{align*}

Note that, for this case, $n$ cannot be 7.

\noindent\textbf{Case (ii). $\alpha_{0,1}=2$.} In this case, we have the following correspondence:
\[(\ldots,3^{\alpha_{0,3}},1^2)\to (\ldots,3^{\alpha_{0,3}},2).\]
Here, in the image, $\sum_{k\geq 0}(2\alpha_{k,6}+\alpha_{k,3})>1$.

\noindent\textbf{Case (iii). $\alpha_{0,1}\geq 3$ and $\alpha_{0,1}=\sum_{k\geq 0}(2\alpha_{k,6}+\alpha_{k,3})-2$.} 

\noindent\textbf{Case (iii-a). $\alpha_{0,3}\geq 2$, or $\sum_{k\geq 1}\alpha_{k,3}\geq 1$.} Here, we have the following correspondence:
\begin{align*}
(\ldots,3^{\alpha_{0,3}},1^{\alpha_{0,1}})\to
\begin{cases}
&(\ldots,5,3^{\alpha_{0,3}-2},2,1^{\alpha_{0,1}-1}),\ \text{if } \alpha_{0,3}\geq 2;\\
&(\ldots,12k^\prime+5,(12k^\prime+3)^{\alpha_{k^\prime,3}-1},\ldots,3^0,2,1^{\alpha_{0,1}-1}),\\
&\ \text{if $\alpha_{0,3}=1$, and $k^\prime$ is the smallest $k\geq1$ for which $\alpha_{k,3}\neq0$}.
\end{cases}
\end{align*}
For example,
\[(\ldots,3^6,1^4)\to (\ldots,5,3^4,2,1^3).\]
Note that for the right-hand side partition, $\alpha_{0,1}=\sum_{k\geq 0}(2\alpha_{k,6}+\alpha_{k,3})-1$.

\noindent\textbf{Case (iii-b). $\alpha_{0,3}=1$, and $\sum_{k\geq 1}\alpha_{k,3}=0$.}
In the above cases, where there is a part of the form $12k+9$ in the image partitions, there is no part of the form $12k+6$. In this case, either there is at least one part of the form $12k+6$ with multiplicity 2, or there are at least two distinct parts of the form $12k+6$ in the pre-image partitions. Therefore, in this case, we add 3 to the smallest part of the form $12k+6$, to get a distinct image. Moreover, we transform $1^{\alpha_{0,1}}$ as follows.
\begin{align*}
	(1^{\alpha_{0,1}})\to\begin{cases}
		\left(3^{\frac{\alpha_{0,1}-3}{3}},2,1\right), &\text{if } \alpha_{0,1}\equiv 0 \pmod 3;\\
		\left(3^{\frac{\alpha_{0,1}-1}{3}},1\right), &\text{if } \alpha_{0,1}\equiv 1 \pmod 3;\\
		\left(3^{\frac{\alpha_{0,1}-5}{3}},2^2,1\right), &\text{if } \alpha_{0,1}\equiv 2 \pmod 3.
	\end{cases}
\end{align*}
For example,
\[(\ldots,6^2,3,1^3)\to (\ldots,9,6,2,1).\]

\noindent\textbf{Case (iv). $\alpha_{0,1}\geq 3$ and $\alpha_{0,1}=\sum_{k\geq 0}(2\alpha_{k,6}+\alpha_{k,3})-1$.} 

\noindent\textbf{Case (iv-a). $\sum_{k\geq 0}\alpha_{k,6}\geq 1$.}
Here, we add 1 to the smallest part of the form $12k+6$, and transform $(1^{\alpha_{0,1}})$ to $(2,1^{\alpha_{0,1}-3})$. For example,
\[(6^2,3,1^4)\to (7,6,3,2,1).\]
Here, for the right-hand side partition, $\alpha_{0,1}=\sum_{k\geq 0}(2\alpha_{k,6}+\alpha_{k,3})-2$. Therefore, the the right-hand side partitions in this case are different from the right-hand side partitions in the other cases.

\noindent\textbf{Case (iv-b). $\sum_{k\geq 0}\alpha_{k,6}=0$.} Here, $\sum_{k\geq 0}\alpha_{k,3}\geq 4$. We add three smallest parts of the form $12k+3$ (all of which can be distinct or equal) to make the sum a part of the form $12k+9$. Moreover, we transform $1^{\alpha_{0,1}}$ similar to \textbf{Case (iii-b)}. For example,
\[(3^4,1^3)\to (9,3,2,1) \text{ and } (3^5,1^4)\to (9,3^3,1).\]
The resulting right-hand side partitions in this case contain no part of the form $12k+6$ and the number of parts of the form $12k+3$ is higher, and therefore they are different from the right-hand side partitions in \textbf{Case (i-a)}, \textbf{Case (i-b)}, and \textbf{Case (iii-b)}.
\section{Concluding Remarks}
Let $t\geq 4$ and $k\geq 0$ be integers. To handle Conjecture \ref{con1} for $t\geq4$, we write the parts of a partition in $\mathcal{B}_t(n)$ modulo $t(t+1)=:\mathbcal{S}_k$. We arrange the general blocks $1_k,~2_k,\ldots,~t_k$ (from below to above) of any partition in $\mathcal{B}_t(n)$ as follows.

\begin{align*}
	\bigg(	& \mathbcal{S}_k+t(t+1)-1, \mathbcal{S}_k+t(t+1)-2, \ldots , \mathbcal{S}_k+t^2+1,  \textcolor{green}{(\mathbcal{S}_k+t^2)^0},\\
	&\textcolor{red}{\big(\mathbcal{S}_k+(t-1)(t+1)\big)}^{\beta_{k,(t-1)(t+1)}},   \ldots ,  \textcolor{green}{(\mathbcal{S}_k+(t-1)t)^0}, \big(\mathbcal{S}_k+(t-1)t-1\big)^{\beta_{k,(t-1)t-1}},\\
	&\hspace{20mm} \vdots \hspace{70mm} \vdots\\
	&\textcolor{red}{\big(\mathbcal{S}_k+i(t+1)\big)}^{\beta_{k,i(t+1)}},   \ldots ,  \textcolor{green}{(\mathbcal{S}_k+it)^0},  \big(\mathbcal{S}_k+it-1\big)^{\beta_{k,it-1}},  \ldots ,  \mathbcal{S}_k+(i-1)t+i,\\
	& \hspace{20mm}\vdots  \hspace{70mm} \vdots\\
	&\textcolor{red}{\big(\mathbcal{S}_k+2(t+1)\big)}^{\beta_{k,2(t+1)}},   \mathbcal{S}_k+2t+1,  \textcolor{green}{(\mathbcal{S}_k+2t)^0}, \ldots , \mathbcal{S}_k+t+2,\\
	&\textcolor{red}{\big(\mathbcal{S}_k+(t+1)\big)}^{\beta_{k,t+1}},   \textcolor{green}{(\mathbcal{S}_k+t)^0}, \big(\mathbcal{S}_k+t-1\big)^{\beta_{k,t-1}},  \ldots ,  \mathbcal{S}_k+2,  \mathbcal{S}_k+1 \bigg)_{k\ge 0},
\end{align*}

\noindent where $\beta_{k,j}$ is the multiplicity of the part $\mathcal{S}_k+j$ and multiplicities of some parts are omitted for simplicity. We can define a map similar to the map $\Phi$ from $\mathcal{B}_t(n)$ to $\mathcal{B}_{t+1}(n)$, call it $\Phi^{*}$. Here, the $(t-1)$ parts marked in red are divisible by $t+1$, which are not supposed to be present in the image partitions. The parts marked in green are divisible by $t$, which are not present in the blocks of a partition from $\mathcal{B}_t(n)$.
\par 
Upon applying the map $\Phi^{*}$, we observe that in some cases the blocks $3_k,\ldots,(t-1)_k$ lose two 2-hooks, when the multiplicity of the part immediately to the right of the red part is zero. It seems possible to find a new map similar to $\Psi$ to compensate these 2-hooks, by which the Conjecture \ref{con1} will be proved for all $t\geq4$.


\end{document}